\documentclass{amsart}
\usepackage{amsmath,amsthm,amssymb,color}
\usepackage[dvipdfm,
 bookmarkstype=toc=true,
 linktocpage=true,
 bookmarks=true
 ]{hyperref}

\theoremstyle{plain} \numberwithin{equation}{section}
\newtheorem{theo}{Theorem}[section]
\def\C{\mathbb C}
\def\R{\mathbb R}
\def\Z{\mathbb Z}
\DeclareMathOperator{\uR}{\underline{\mathbb{R}}}
\DeclareMathOperator{\uC}{\underline{\mathbb{C}}}
\DeclareMathOperator{\ue}{\underline{{e}}}

\DeclareMathOperator{\GL}{GL}
\begin{document}
	\title{Invariant stably complex structures on topological toric manifolds}
	\date{\today}
	\author{Hiroaki Ishida}
	\address{Department of Mathematics, Graduate School of Science, Osaka City University, Sugimoto, Sumiyoshi-ku, Osaka 558-8585, Japan}
	\email{hiroaki.ishida86@gmail.com}
	\begin{abstract}
		We show that any $(\C ^*)^n$-invariant stably complex structure on a topological toric manifold of dimension $2n$ is integrable. 
		We also show that such a manifold is weakly $(\C ^*)^n$-equivariantly isomorphic to a toric manifold.  
	\end{abstract}
	\maketitle
	\section{Introduction}
		A \emph{toric manifold} is a nonsingular complete toric variety. As a topological analogue of a toric manifold, the notion of \emph{topological toric manifold} has been introduced by the author, Y. Fukukawa and M. Masuda \cite{IFM}. 
		A topological toric manifold of dimension $2n$ is a smooth closed manifold endowed with an effective $(\C ^*)^n$-action having an open dense orbit, and  locally equivariantly diffeomorphic to a \emph{smooth} representation space of $(\C ^*)^n$. 
We note that a topological toric manifold is locally equivariantly diffeomorphic to an \emph{algebraic} representation space if and only if it is a toric manifold. 

		A \emph{quasitoric manifold} introduced by M. Davis and T. Januskiewicz \cite{DJ91} of dimension $2n$ is a smooth closed manifold endowed with a locally standard $(S^1)^n$-action, whose orbit space is a simple polytope. In \cite{IFM}, it is shown that any quasitoric manifold is a topological toric manifold with the restricted compact torus action.  Conversely, it is also shown that any topological toric manifold of dimension less than or equal to $6$ with the restricted compact torus action is a quasitoric manifold. However, there are infinitely many topological toric manifolds with the restricted compact torus action which are not equivariantly diffeomorphic to any quasitoric manifold.

	Among quasitoric manifolds, some admit invariant almost complex structures under the compact torus actions. M. Masuda provided examples of $4$-dimensional quasitoric manifolds which admit $(S^1)^2$-invariant almost complex structures (see \cite[Theorem 5.1]{Mas99}). A. Kustarev described a necessary and sufficient condition for a quasitoric manifold to admit a torus invariant almost complex structure for arbitrary dimension (see \cite[Theorem 1]{Kus09}). 
	 
	As we mentioned, any quasitoric manifold is a topological toric manifold with the restricted compact torus action.  In this paper, we discuss on $(\C ^*)^n$-invariant stably, or almost complex structures on topological toric manifolds of dimension $2n$. The followings are our results: 
	\begin{theo}\label{MainTheo}
		Let $X$ be a topological toric manifold of dimension $2n$. Let $\uR ^{2\ell}$ be the product bundle of rank $2\ell$ over $X$, $TX$ the tangent bundle of $X$. 
If there exists a $(\C ^*)^n$-invariant stably complex structure $J$ on $TX \oplus \uR ^{2\ell}$, then $TX$ becomes a complex subbundle of $TX \oplus \uR ^{2\ell}$. Namely, $X$ has an invariant almost complex structure. 
	\end{theo}
	\begin{theo}\label{MainTheo2}
		Let $X$ be a topological toric manifold of dimension $2n$, $J$ a $(\C ^*)^n$-invariant almost complex structure. Then, $J$ is integrable and $X$ is weakly equivariantly isomorphic to a toric manifold. Namely, there are a toric manifold $Y$, a biholomorphism $f: X \to Y$ and a smooth automorphism $\rho$ of $(\C ^*)^n$ such that $f\circ g = \rho (g)\cdot f$ for all $g \in (\C ^*)^n$. 
	\end{theo}
	
	If we replace the condition \lq\lq{$(\C ^*)^n$-invariant}\rq\rq by \lq\lq{$(S^1)^n$-invariant}\rq\rq on the almost complex structure $J$, then Theorem \ref{MainTheo2} does not hold. For example, $\C P^2 \# \C P^2 \# \C P^2$ with an effective $(S^1)^2$-action is a topological toric manifold with the restricted $(S^1)^2$-action. One can show that there  exists an $(S^1)^2$-invariant almost complex structure on $\C P^2 \# \C P^2 \# \C P^2$ (see \cite[Theorem 5.1]{Mas99}). However, $\C P^2 \# \C P^2 \# \C P^2$ carries no complex structure because $\C P^2 \# \C P^2 \# \C P^2$ does not fit Kodaira's classification of complex surfaces. Namely, the almost complex structure is not integrable.  

	For a topological toric manifold $X$ of dimension $2n$, there is a canonical short exact sequence of $(\C ^*)^n$-bundles 
	\begin{equation*}
		0 \to \uC ^{m-n} \to \bigoplus _{i =1}^m L_i \to TX \to 0,
	\end{equation*}
	where 
	$L_i$'s are complex line bundles.  
	(see \cite[Theorem 6.1]{IFM}). Theorems \ref{MainTheo} and \ref{MainTheo2} say 
	that the short exact sequence above does not split as $(\C ^*)^n$-bundles unless $X$ is a toric manifold.
	\section{Preliminaries}\label{sec:preliminary}
		In this section, we review the quotient construction of topological toric manifolds and the correspondence between topological toric manifolds and nonsingular complete topological fans (see \cite{IFM} for details). 
		
		A \emph{nonsingular complete topological fan} is a pair $\Delta = (\Sigma ,\beta )$ such that
		\begin{enumerate}
			\item $\Sigma$ is an abstract simplicial complex on $[m]=\{ 1,\dots ,m\}$,
			\item $\beta : [m] \to (\C \times \Z )^n$ is a function which satisfies the following: 
			\begin{enumerate}
				\item \label{compfan}Let $\mathrm{Re}$ be the composition of two natural projections $(\C \times \Z)^n \to \C ^n$ and $\C ^n \to \R ^n$. We assign a cone 
					\begin{equation*}
						\left\{ \sum _{i \in I}a_i(\mathrm{Re}\circ \beta )(i) \mid a_i \geq 0\right\} 
					\end{equation*}to each simplex $I$ in $\Sigma $. 
				Then, we have a collection of cones in $\R ^n$. 
				\begin{enumerate}
					\item Each pair of two cones does not overlap on their relative interiors. Namely, the real part $\mathrm{Re} \circ \beta$ of $\beta$ together with $\Sigma$ forms an ordinary fan.
					\item The union of all cones coincides with $\R ^n$. Namely, the fan is complete. 
				\end{enumerate}
				\item The integer part of $\beta$ together with $\Sigma$ forms a nonsingular multi-fan (see \cite[p.249]{Mas99}).
			\end{enumerate}
		\end{enumerate}
		It follows from (\ref{compfan}) that $\Sigma $ must be a simplicial $(n-1)$-sphere with $m$ vertices. If we regard integers $\Z$ as a subset of $\C \times \Z$ via $a \mapsto (a,a)$ for $a \in \Z$, then any nonsingular complete fan can be regarded as a special case of a nonsingular complete topological fan. Conversely, if the image of $\beta $ is contained in the diagonal subgroup $\Z ^n$ of $(\C \times \Z )^n$, then $\Delta $ becomes a nonsingular complete fan.    
		
		We express $\beta (i)$ as $\beta _i = (\beta _i^1,\dots ,\beta _i^n) \in (\C \times \Z )^n$ and $\beta _i^j = (b_i^j+\sqrt{-1}c_i^j,v_i^j) \in \C \times \Z$ for $i= 1,\dots ,m$ and $j=1,\dots ,n$. For a nonsingular complete topological fan $\Delta = (\Sigma , \beta )$, we can construct a topological toric manifold as follows. We set 
			\begin{equation*}
				U(I) := \{ (z_1,\dots ,z_m) \in \C ^m \mid z_i \neq 0 \text{ for $i \notin I$} \}
			\end{equation*}
			for $I \in [m]$, and 
			\begin{equation*}
				U(\Sigma ) := \bigcup _{I \in \Sigma }U(I).
			\end{equation*}
			We define a group homomorphism $\lambda _{\beta _i} : \C ^* \to (\C ^*)^n$ by 
			\begin{equation}\label{eq:lambdabeta}
				\lambda _{\beta _i}(h_i) := (h_i^{\beta _i^1},\dots ,h_i^{\beta _i^n}),
			\end{equation}
			where 
			\begin{equation}\label{eq:betapower}
				h_i^{\beta _i^j} := |h_i|^{b_i^j+\sqrt{-1}c_i^j}\left( \frac{h_i}{|h_i|}\right) ^{v_i^j}.  
			\end{equation}
			We define a group homomorphism $\lambda : (\C ^*)^m \to (\C ^*)^n$ by 
			\begin{equation*}
				\lambda (h_1,\dots ,h_m) := \prod _{i =1}^m \lambda _{\beta _i}(h_i).
			\end{equation*}
			One can show that $\lambda $ is surjective (see \cite[Lemma 4.1]{IFM}). 
We note that the $(\C ^*)^m$-action on $U(\Sigma )$ given by coordinatewise multiplications induces the action of $(\C ^*)^m/\ker \lambda $ on the quotient space $X(\Delta ) := U(\Sigma )/\ker \lambda $. Since $\lambda $ is surjective, we can identify $(\C ^*)^m/\ker \lambda$ with $(\C ^*)^n$ through $\lambda $. Hence $X(\Delta )$ is equipped with the $(\C ^*)^n$-action. One can show that $X(\Delta )$ is a topological toric manifold. 
		
		We shall remember the equivariant charts and transition functions of $X(\Delta )$ described in \cite{IFM} for later use.
		We regard $\beta _i^j = (b_i^j+\sqrt{-1}c_i^j, v_i^j)$ as the following matrix: 
		\begin{equation*}
			\beta _i^j =
			\begin{pmatrix}
				b_i^j & 0\\
				c_i^j & v_i^j
			\end{pmatrix}.
		\end{equation*}
		And we also regard $\beta _i$ as an $n$-tuple $(\beta _i^1,\dots ,\beta _i^n)$ of square matrices of size $2$.  Let $\Sigma ^{(n)}$ denote the set of $(n-1)$-dimensional simplices in $\Sigma$. For $I \in \Sigma ^{(n)}$, the dual $\{ \alpha _i^I\} _{i \in I}$ of $\{ \beta _i\} _{i \in I}$ is defined to be 
		\begin{equation}\label{eq:alphabeta}
			\langle \alpha _h^I, \beta _i \rangle = \delta _h^i \mathbf{1}, 
		\end{equation}
		where $\delta$ denotes the the Kronecker delta, and $\langle \ , \  \rangle$ is given by 
		\begin{equation*}
			\langle \alpha , \beta \rangle = \sum _{j=1}^n \alpha ^j \beta ^j
		\end{equation*}
		for $n$-tuples $\alpha = (\alpha ^1,\dots ,\alpha ^n)$ and $\beta = (\beta ^1,\dots ,\beta ^n)$ of square matrices of size $2$. The equivariant charts are given as follows. 
		For $\alpha = (\alpha ^1,\dots ,\alpha ^n)$, we define a representation $\chi ^{\alpha } : (\C ^*)^n \to \C ^*$ by 
		\begin{equation*}
			\chi ^{\alpha} (g_1,\dots ,g_n) := \prod _{j =1}^n g_j^{\alpha ^j}.
		\end{equation*}
		Let $V(\chi ^{\alpha})$ denote the representation space of $\chi ^{\alpha}$. For $I \in \Sigma ^{(n)}$, the equivariant chart $\varphi _I : U(I)/\ker \lambda \to \bigoplus _{i \in I}V(\chi ^{\alpha _i^I})$ is defined by 
		\begin{equation*}
			\varphi _I([z_1,\dots ,z_m]) := \left (\prod _{j=1}^m z_j ^{\langle \alpha _i^I,\beta _j\rangle } \right) _{i \in I}, 
		\end{equation*}
		where $[z_1,\dots , z_m]$ denotes the equivalence class of $(z_1,\dots ,z_m) \in U(\Sigma )$. 
		Let $(w_i)_{i \in I}$ be the coordinates of $\bigoplus _{i \in I}V(\chi ^{\alpha _i^I})$. Namely, 
		\begin{equation*}
			w_i := \prod _{j=1}^m z_j ^{\langle \alpha _i^I,\beta _j\rangle }.
		\end{equation*}

		An \emph{omniorientation} of a topological toric manifold $X$ is a choice of orientations of normal bundles of \emph{characteristic} submanifolds of $X$. 
Here, a characteristic submanifold of $X$ is a connected $(\C ^*)^n$-invariant submanifold of codimension $2$. By the construction above, $\beta $ allows us to decide an omniorientation of $X(\Delta )$ as follows. Each characteristic submanifold intersecting $U(I)/\ker \lambda $ is locally represented as $\{ (w_i)_{i \in I} \mid w_{i_0} = 0\}$ for some $i_0 \in I$ through $\varphi _I : U(I)/\ker \lambda \to \bigoplus _{i \in I}V(\chi ^{\alpha _i^I})$. We choose the orientation of the normal bundle of each  characteristic submanifold so that $\varphi _I$ preserves the orientations of the normal bundles of the characteristic submanifolds and 
the normal bundle of $\{ (w_i)_{i \in I} \mid w_{i_0} = 0\}$ in $\bigoplus _{i \in I}V(\chi ^{\alpha _i^I})$. 
We need to show that this is well-defined. By \eqref{eq:alphabeta}, $\lambda _{\beta _{i_0}}(\C ^*)$ fixes $\{ (w_i)_{i \in I} \mid w_{i_0} = 0\}$. Let us consider the differential $(\lambda _{\beta _{i_0}} (\sqrt{-1})) _*$ of the action $\lambda _{\beta _{i_0}} (\sqrt{-1})$ at a point of $\{ (w_i)_{i \in I} \mid w_{i_0} = 0\}$. The normal vector space of $\{ (w_i)_{i \in I} \mid w_{i_0} = 0\}$  in $\bigoplus _{i \in I}V(\chi ^{\alpha _i^I})$  at any point  in $\{ (w_i)_{i \in I} \mid w_{i_0} = 0\}$  can be identified with $V(\chi ^{\alpha _{i_0}^I})$. By direct computation, we have 
		\begin{equation*}
			\chi ^{\alpha _{i_0}^I} \circ \lambda _{\beta _{i_0}}(\sqrt{-1}) = \sqrt{-1}
		\end{equation*}
		(see \cite[Lemma 2.2.(2)]{IFM}). Namely, for any point of the characteristic submanifold and any nonzero normal vector $\xi $ at the point, if we choose the orientation of the normal bundle of the characteristic submanifold so that ($\xi, (\lambda _{\beta _{i_0}} (\sqrt{-1})) _*(\xi))$ is a positive basis, then $\varphi _I$ preserves the orientations of the normal bundles. 
		
The correspondence $\Delta \mapsto X(\Delta )$ is bijective between nonsingular complete topological fans and omnioriented topological toric manifolds (see \cite[Theorem 8.1]{IFM}). 
		We shall see the inverse correspondence. For a topological toric manifold $X$ of dimension $2n$ with an omniorientation, let us denote characteristic submanifolds of $X$ by $X_1,\dots ,X_m$. Define 
		\begin{equation*}
			\Sigma = \left\{ I \in [m] \mid \bigcap _{i \in I} X_i \neq \emptyset \right\}.
		\end{equation*}
		For an orientation on normal bundle of $X_i$, we can find a unique complex structure $J_i$ such that 
		\begin{itemize}
			\item the orientation coincides with the orientation which comes from $J_i$, 
			\item the $\C ^*$-subgroup of $(\C ^*)^n$ which fixes each point of $X_i$ acts on the normal bundle as $\C$-linear with respect to $J_i$.   
		\end{itemize}
		For  $J_i$, we can find a unique $\beta _i \in (\C \times \Z)^n$ such that 
		\begin{equation*}
				(\lambda _{\beta _i}(h))_*(\xi )= h \xi  
		\end{equation*}
		for any normal vector $\xi $ and $h \in \C ^*$, where the right hand side is the multiplication with complex number. For an omniorientation of $X$, define $\beta : [m] \to (\C \times \Z )^n$ as $\beta (i) := \beta _i$. Then, the pair $\Delta (X) =(\Sigma , \beta )$ becomes a nonsingular complete topological fan and the correspondence $X \mapsto \Delta (X)$ is the inverse correspondence of $\Delta \mapsto X(\Delta )$. Namely, there exists an equivariant diffeomorphism $X \to X(\Delta (X))$ which preserves the omniorientations. 

		The transition functions of $X(\Delta )$ are given as follows. Let $K$ be another element in $\Sigma ^{(n)}$. By direct computation, $k$-component of $\varphi _K(\varphi _I ^{-1}(w_i)_{i \in I})$ for $k \in K$ is given as
		\begin{equation}\label{eq:transition}
			\prod _{i \in I}w_i^{\langle \alpha ^K_k,\beta _i\rangle}
		\end{equation}
		(see \cite[Lemma 5.2]{IFM}). We remark that
		\begin{equation*}
			\frac{\partial}{\partial \bar{w_j}}\left( \prod _{i \in I}w_i^{\langle \alpha ^K_k,\beta _i\rangle}\right) =0
		\end{equation*}
		if and only if $\langle \alpha ^K_k,\beta _j\rangle = \gamma ^K_{k,j}\mathbf{1}$ for some integer $\gamma ^K_{k,j}$ (see \eqref{eq:betapower}). This implies that all transition functions are holomorphic if and only if there is an integer $\gamma ^K_{k,j}$ such that $\langle \alpha ^K_k,\beta _i\rangle = \gamma ^K_{k,j}\mathbf{1}$ for all $i \in [m]$, $k \in K$ and $K \in \Sigma ^{(n)}$. In this case, each transition function is a Laurent monomial and hence $X(\Delta )$ is weakly equivariantly diffeomorphic to a toric manifold.
	\section{Proof of Theorem \ref{MainTheo}}\label{Sec:proof}
		Let $X$ be a $2n$-dimensional topological toric manifold, $TX$ the tangent bundle of $X$, $J$ a $(\C ^*)^n$-invariant complex structure on $TX \oplus \uR ^{2\ell}$. We take an omniorientation of $X$ and consider the topological fan $\Delta = (\Sigma , \beta )$ associated to $X = X(\Delta )$ with the given omniorientation.  We define a cross section $\ue _h: X \to \uR ^{2\ell}=X \times \R ^{2\ell}$ for $h =1,\dots ,2\ell$ by $x \mapsto (x, e_h)$ for all $x \in X$, where $e_h$ denotes the $h$-th standard basis vector of $\R ^{2\ell}$.  

		
		There is a natural inclusion $(\C ^*)^n \hookrightarrow X$ given by $g \mapsto g\cdot [1,\dots ,1]$ where $[1,\dots ,1]$ denotes the equivalence class of $(1,\dots ,1)$ in $U(\Sigma )$. For $I \in \Sigma ^{(n)}$, the inclusion is of the form  
		\begin{equation}\label{eq:inclusion}
			\bigoplus _{i \in I}\chi ^{\alpha _i^I} : g=(g_j)_{j=1,\dots ,n} \mapsto \left( \chi ^{\alpha _i^I}(g) \right) _{i \in I}
		\end{equation}
		via the equivariant local chart $\varphi _I : U_I/\ker \lambda \to \bigoplus _{i \in I}V(\chi ^{\alpha _i^I})$. We identify $\bigoplus _{i \in I}V(\chi ^{\alpha _i^I})$ with $\R ^{2n}$ by 
		\begin{equation}\label{eq:cpxreal}
			w_i = x_i +\sqrt{-1}y_i
		\end{equation}
		for $i \in I$, where $(w_i)_{i \in I}$ denote the coordinates of $\bigoplus _{i \in I} V(\chi ^{\alpha _i^I})$. We also identify $(\C ^*)^n$ with $(\R \times \R /2\pi \Z )^n$ by 
		\begin{equation}\label{eq:pd}
			\psi : (g_j)_{j =1,\dots , n} \mapsto \left(\log |g_j| , -\sqrt{-1}\log \left( \frac{g_j}{|g_j|} \right) \right) _{j=1,\dots ,n} . 
		\end{equation}		
		Let $(\tau _j, \theta _j)_{j =1,\dots ,n}$ be the coordinates of $(\R \times \R /2\pi \Z )^n$. Since $J$ is $(\C ^*)^n$-invariant, the matrix representation, denoted $J_0$, of $J$ on $(\R \times \R /2\pi \Z )^n$  with respect to the coordinates $(\tau _j, \theta _j)_{j =1,\dots ,n}$ and sections $\ue _h$'s is constant. 
		
		Let $\Psi _I : (\R ^2 \setminus \{ 0\} )^n \to (\R \times \R /2\pi \Z )^n$ be the composition of the identification \eqref{eq:cpxreal}, the inverse of \eqref{eq:inclusion}, and $\psi$. Namely, 
		\begin{equation}\label{eq:Psi_I}
			\Psi _I((x_i,y_i)_{i \in I} ) := \psi \circ \left( \bigoplus _{i \in I}\chi ^{\alpha _i^I}\right) ^{-1}((x_i + \sqrt{-1} y_i)_{i \in I}).
		\end{equation}
		Since $(\bigoplus _{i \in I}\chi ^{\alpha _i^I})^{-1}$ coincides with $\prod _{i \in I} \lambda _{\beta _i}$ (see \cite[Lemma 2.3]{IFM}), 
		it follows from \eqref{eq:lambdabeta} and \eqref{eq:betapower} that the coordinates $(\tau _j, \theta _j)_{j=1,\dots ,n}$ are represented as  
		\begin{align*}
			\tau _j &= \log \left( \prod _{i \in I} |(x_i +\sqrt{-1}y_i)^{\beta _i^j}| \right)\\ &=\frac{1}{2}\sum _{i \in I}b_i^j\log (x_i^2+y_i^2)
		\end{align*}
		and 
		\begin{align*}
			\theta _j &= -\sqrt{-1} \log \left( \prod _{i \in I}\frac{(x_i +\sqrt{-1}y_i)^{\beta _i^j}}{|(x_i +\sqrt{-1}y_i)^{\beta _i^j}|}\right) \\ 
			&= -\sqrt{-1}\log \left( \prod _{i \in I}|x_i+\sqrt{-1}y_i|^{\sqrt{-1}c_i^j}\left( \frac{x_i + \sqrt{-1}y_i}{|x_i+ \sqrt{-1}y_i|}\right) ^{v_i^j} \right) \\
			&= \sum _{i \in I}\left( \frac{c_i^j+\sqrt{-1}v_i^j}{2}\log (x_i^2+y_i^2) -\sqrt{-1}v_i^j\log (x_i+\sqrt{-1}y_i) \right).
		\end{align*}
		Then, by direct computation, we have
		\begin{equation*}
			\frac{\partial \tau _j}{\partial x_i} = \frac{b_i^jx_i}{x_i^2+y_i^2}, \quad \frac{\partial \tau _j}{\partial y_i} = \frac{b_i^jy_i}{x_i^2+y_i^2},
		\end{equation*}
		\begin{align*}
			\frac{\partial \theta _j}{\partial x_i} &= \frac{(c_i^j+\sqrt{-1}v_i^j)x_i}{x_i^2+y_i^2}-\sqrt{-1}v_i^j\frac{1}{x_i+\sqrt{-1}y_i} \\
			&= \frac{c_i^jx_i - v_i^jy_i}{x_i^2+y_i^2}
		\end{align*}
		and
		\begin{align*}
			\frac{\partial \theta _j}{\partial y_i} &= \frac{(c_i^j+\sqrt{-1}v_i^j)y_i}{x_i^2+y_i^2}+v_i^j\frac{1}{x_i+\sqrt{-1}y_i} \\
			&= \frac{c_i^jy_i+v_i^jx_i}{x_i^2+y_i^2}.
		\end{align*}
		Therefore, 
		\begin{equation*}
			\begin{pmatrix}
				\frac{\partial \tau _j}{\partial x_i} & \frac{\partial \tau _j}{\partial y_i}\\
				\frac{\partial \theta _j}{\partial x_i} & \frac{\partial \theta _j}{\partial y_i}
			\end{pmatrix}
			=
			\begin{pmatrix}
				b_i^j & 0\\
				c_i^j & v_i^j
			\end{pmatrix}
			\begin{pmatrix}
				\frac{x_i}{x_i^2+y_i^2} & \frac{y_i}{x_i^2+y_i^2}\\
				\frac{-y_i}{x_i^2+y_i^2} & \frac{x_i}{x_i^2+y_i^2}
			\end{pmatrix}
			=\beta _i^jt_i,
		\end{equation*}
		where 
		\begin{equation*}
			t_i = \frac{1}{x_i^2+y_i^2}\begin{pmatrix}
				x_i & y_i \\
				-y_i & x_i
			\end{pmatrix} \in \C ^* \subset \GL (2,\R).
		\end{equation*}
		
		We set two square matrices 
		\begin{equation*}
			B=(\beta _i^j)_{i,j} \text{ and } T=\mathrm{diag}(t_i ;i\in I)
		\end{equation*}
		of size $n$ whose entries are square matrices of size $2$. 
		Then, the differential $T_{(x,y)}\Psi _I$ of $\Psi _I$ at $(x,y)$ is represented as $BT$ with respect to the coordinates $(x_i,y_i)_{i \in I}$ and $(\tau _j,\theta _j)_{j=1,\dots ,n}$.
		Hence the complex structure $J$ of $TX \oplus \uR ^{2\ell}$ is represented on $\bigoplus_{i\in I}V(\chi^{\alpha_i^I})=\R^{2n}$ as the following square matrix 
		\begin{equation}\label{eq:J_I}
			J_I := \begin{pmatrix}
				(BT)^{-1} &  \\
				& I_{2\ell}
			\end{pmatrix}J_0\begin{pmatrix}
				BT &  \\
				& I_{2\ell}
			\end{pmatrix}
		\end{equation}
		of size $2n + 2\ell$ with respect to the coordinates $(x_i,y_i)_{i \in I}$ and sections $\ue _h$, where $I_{2\ell}$ denote the identity matrix of size $2\ell$.
		We set 
		\begin{equation*}
			J_0 =: \begin{pmatrix}
				J_{11} & J_{12}\\
				J_{21} & J_{22}
			\end{pmatrix}
		\end{equation*}
		where $J_{11}, J_{12}, J_{21}, J_{22}$ are matrices of $2n\times 2n, 2n \times 2\ell , 2\ell \times 2n , 2\ell \times 2\ell$, respectively. Then, 
		\begin{equation}\label{eq:J_I'}
			J_I = \begin{pmatrix}
				T^{-1}(B^{-1}J_{11} B)T & T^{-1}B^{-1}J_{12}\\
				J_{21}BT & J_{22}
			\end{pmatrix}.
		\end{equation}
		Each entry of $J_I$ is a homogeneous (rational) polynomial in $(x_i,y_i)_{i \in I}$. However, each entry of $J_I$ must be a smooth function on $\R ^{2n}$, in particular, at the origin. 
		Hence each entry of $J_{21}$ must be $0$. This means that the tangent space at any point of $X$ is stable under $J$. 
		It follows that $TX$ is a complex subbundle of $TX \oplus \uR ^{2\ell}$ with respect to $J$. The theorem is proved. \qed
	\section{Proof of Theorem \ref{MainTheo2}}
		Let $X$ be a topological toric manifold of dimension $2n$ with a $(\C ^*)^n$-invariant almost complex structure $J$. Since any characteristic submanifold of $X$ and $X$ itself are almost complex submanifolds, the normal bundles of characteristic submanifolds of $X$ become complex line bundles. Hence, we have a topological fan $\Delta =(\Sigma ,\beta )$ associated to $X$. Namely, for each characteristic submanifold $X_i$ of X, we choose the unique $\beta _i \in (\C \times \Z )^n$ so that 
		$\lambda _{\beta _i}(\C ^*)$ fixes all points in $X_i$, and 
		\begin{equation*}
			(\lambda _{\beta _i}(h))_*(\xi )=h\xi 
		\end{equation*}
		for any $h \in \C ^*$ and any normal vector $\xi $ of $X_i$. 
		
		According to the proof of Theorem \ref{MainTheo}, we identify the dense orbit with $(\R \times \R /2\pi \Z )^n$ and $\bigoplus _{i \in I}V (\chi ^{\alpha _i^I})$ with $\R ^{2n}$ for $I \in \Sigma ^{(n)}$. Since $J$ is $(\C ^*)^n$-invariant, the matrix representation, denoted $J_0$, of $J$ on $(\R \times \R /2\pi \Z )^n$ with respect to the coordinate $(\tau _j, \theta _j)_{j=1,\dots ,n}$ of $(\R \times \R /2\pi \Z )^n$ is constant. Let us remind $\Psi _I : (\R ^2 \setminus \{ 0\})^n \to (\R \times \R /2\pi \Z)^n$ (see \eqref{eq:Psi_I}). Then, the almost complex structure $J$ is represented on $\bigoplus _{i \in I}V (\alpha ^I_i) = \R ^{2n}$ as the following square matrix
		\begin{equation}\label{eq:J_I''}
			J_I := (BT)^{-1}J_0(BT) = T^{-1}(B^{-1}J_0B)T
		\end{equation}
		(this is the case when $\ell =0$ in \eqref{eq:J_I} and \eqref{eq:J_I'}).
		
		As we saw at the end of Section \ref{Sec:proof}, each entry of $J_I$ is a homogeneous (rational) polynomial in $(x_i,y_i)_{i \in I}$. However, each entry of $J_I$ must be a smooth function on $\R ^{2n}$, in particular, at the origin. Hence $J_I$ must be constant. Since $J_I$ and $B^{-1}J_0B$ are both constant and $T$ can take any element in $(\C ^*)^n$, it follows from \eqref{eq:J_I''} that $J_I$ and $B^{-1}J_0B$ must be both in $(\C ^*)^n$ and hence $J_I = B^{-1}J_0B$. Moreover, $J_I^2$ is the minus identity matrix because $J_I$ is the matrix representing the almost complex structure $J$.
		It follows that $J_I$ must be of the form
%
		\begin{equation*}
			\begin{pmatrix}
				\pm \begin{pmatrix}
					0 & -1\\
					1 & 0
				\end{pmatrix} & & \\
				& \ddots & \\
				& & \pm \begin{pmatrix}
					0 & -1\\
					1 & 0
				\end{pmatrix}
			\end{pmatrix}.
		\end{equation*}
		However, each $\beta _i$ is taken so that $(\lambda _{\beta _i}(h))_*(\xi )=h\xi $ for any $h \in \C ^*$ and any normal vector $\xi $ of $X_i$, as we mentioned in Section \ref{sec:preliminary}, we have
		\begin{equation*}
			J_I = \begin{pmatrix}
				\begin{pmatrix}
					0 & -1\\
					1 & 0
				\end{pmatrix} & & \\
				& \ddots & \\
				& & \begin{pmatrix}
					0 & -1\\
					1 & 0
				\end{pmatrix}
			\end{pmatrix}.
		\end{equation*}
		Clearly, the complex structure $J_I$ on $\R ^{2n}$ comes from the identification \eqref{eq:cpxreal}. Therefore, $J$ is integrable and the local chart $\varphi _I: U_I/\ker \lambda \to \bigoplus _{i \in I}V(\chi ^{\alpha_i^I})$ is a holomorphic chart. This implies that for another simplex $K \in \Sigma ^{(n)}$, $k$-component of $\varphi _K(\varphi _I ^{-1}(w_i)_{i \in I})$ given as \eqref{eq:transition}  for $k \in K$ is holomorphic. Thus, the transition functions must be Laurent monomials as remarked at the end of Section \ref{sec:preliminary} and hence $X(\Delta )$ is weakly equivariantly isomorphic to a toric manifold. The theorem is proved. \qed

		\bigskip
		{\bf Acknowledgement.} The author would like to thank Professor Mikiya Masuda for stimulating discussion and helpful comments on the presentation of the paper.  
	

\begin{thebibliography}{99}
		\bibitem{DJ91}
			M. Davis and T. Januskiewicz, \textit{Convex polytopes, Coxeter orbifolds and torus actions}, Duke Math. J. 62 (1991), 417--451. 
		\bibitem{Kus09}
			A. Kustarev, \textit{Equivariant almost complex structures on quasitoric manifolds}, Russian Math. Surveys 64 (2009), no. 1, 156--158. 
		\bibitem{IFM}
			H. Ishida, Y. Fukukawa and M. Masuda, \textit{Topological toric manifolds}, preprint, arXiv:math.AT/10121786.
		\bibitem{Mas99}
			M. Masuda, \textit{Unitary toric manifolds, multi-fans and equivariant index}, Tohoku Math. J. 51 (1999), 237--265.
	\end{thebibliography}
\end{document}